\documentclass[a4paper,12pt,reqno]{amsart}

\usepackage{amssymb,amsmath,amsthm,amsfonts}
\usepackage{graphicx}
\usepackage{fullpage}

\theoremstyle{plain}
\newtheorem{theorem}{Theorem}[section]
\newtheorem*{theorem*}{Theorem}
\newtheorem{proposition}[theorem]{Proposition}
\newtheorem*{proposition*}{Proposition}
\newtheorem{corollary}[theorem]{Corollary}
\newtheorem*{corollary*}{Corollary}
\newtheorem{lemma}[theorem]{Lemma}
\newtheorem*{lemma*}{Lemma}

\theoremstyle{definition}
\newtheorem{definition}[theorem]{Definition}
\newtheorem*{definition*}{Definition}

\newtheorem*{definitions*}{Definitions}

\newtheorem*{example*}{Example}

\newtheorem*{examples*}{Examples}

\newtheorem*{exercise*}{Exercise}

\newtheorem*{exercises*}{Exercises}
\newtheorem{hypothesis}[theorem]{Hypothesis}
\newtheorem*{hypothesis*}{Hypothesis}

\theoremstyle{remark}

\newtheorem*{remark*}{Remark}

\newtheorem*{remarks*}{Remarks}

\newtheorem*{note*}{Note}

\newtheorem*{notes*}{Notes}

\newtheorem*{notation*}{Notation}

\numberwithin{equation}{section}

 \def\N{\mathbb N}  \def\C{\mathbb
  C} \def\R{\mathbb R}  
\def\M{\mathbb M} \def\fP{\mathfrak P} \def\fQ{\mathfrak Q}
\def\eps{\varepsilon} \def\cC{\mathcal C}

\title{Some questions concerning attractors for non-autonomous
  dynamical systems}

\author[R. Johnson]{Russell Johnson}

\address[R. Johnson]{Dipartimento di Sistemi e Informatica, Facolt\`a
  di Ingegneria, Universit\`a di Firenze, Via S. Marta 3, 50139
  Firenze, Italy.}

\email[Russell Johnson]{johnson@dsi.unifi.it}

\author[V. Mu\~{n}oz-Villarragut]{V\'{\i}ctor Mu\~{n}oz-Villarragut}

\address[V. Mu{\~n}oz-Villarragut]{Departamento de Matem\'{a}tica
  Aplicada, E.T.S. de Ingenieros Industriales, Universidad de
  Valladolid, Paseo del Cauce s/n, 47011 Valladolid, Spain.}

\email[V\'{\i}ctor Mu\~{n}oz-Villarragut]{vicmun@wmatem.eis.uva.es}

\date{}

\keywords{Non-autonomous dynamical system, attractor, Pliss reduction principle}

\subjclass{37B55, 34C29, 34C45, 34D45}

\begin{document}

\vspace*{-1.5cm}

\noindent\begin{tabular}{|l|}
  \hline
  Russell Johnson, Víctor Muñoz-Villarragut,\\
  Some questions concerning attractors for non-autonomous dynamical systems,\\
  Nonlinear Analysis: Theory, Methods \& Applications {\bfseries 71} (12), 2009, e1858-e1868.\\
  https://doi.org/10.1016/j.na.2009.02.088\\
  \copyright ~Elsevier\\
  \hline
\end{tabular}

\vspace{1cm}

\begin{abstract}
We compare various concepts of attractor in the context of
non-autonomous dynamical systems. Then, we prove an appropriate
version of the Pliss reduction principle for non-autonomous
differential systems with rapidly oscillating coefficients.
\end{abstract}

\maketitle

\section{Introduction}

An important theme in the theory of non-autonomous differential and
discrete systems is that of integral manifolds and the invariant sets
which they contain. Questions concerning the existence and properties
of attracting invariant sets are of interest. In particular, it is
important to know when an invariant set which is locally attracting
with respect to the restriction of a differential system to the
integral manifold, is actually locally attracting with respect to the
full differential system. That is, one wants to know if the ``Pliss
reduction principle'' (see~\cite{plis}) is valid.\par

In this paper, we consider two topics concerning integral manifolds
for non-autonomous differential systems and the attractors which they
may contain. In our discussion, we will adopt the framework of
skew-product flows (\cite{bebu},\cite{sell}), which has been found to
be convenient in the study of a wide spectrum of problems concerning
non-autonomous differential systems.\par

The first topic is that of the relationship between diverse notions of
``attractor'' in the skew-product framework. We adopt a Lyapunov type
definition of attractor as starting point, and compare it with a
definition of \emph{pullback attractor} which seems appropriate for
skew-product flows. This last concept has been widely discussed and
applied in the recent literature on non-autonomous dynamical systems
(see e.g. \cite{carc}, \cite{sche}, \cite{schm}). Still another type
of attractor is determined in certain circumstances as a fixed point
of a contraction operator defined using a given skew-product flow; see
e.g. \cite{fajm}, \cite{joma} for examples of such ``fixed-point''
attractors. We will see that, under mild hypotheses, an attractor is
of Lyapunov type if and only if it is of pullback type. We will also
see that a fixed-point attractor is of pullback type, and so of
Lyapunov type as well.\par

The second issue we address concerns the behavior of solutions of
non-autonomous differential systems with \emph{rapidly varying
  coefficients}, which lie in or near an integral manifold. The
construction of integral manifolds for such systems was studied by
Coppel and Palmer (\cite{copa},\cite{palm}). More recently this theme
has been taken up by Cheban, Duan and Gherco~\cite{chdg}, and by
Fabbri, Johnson and Palmer~\cite{fajp}. Our point of departure is the
theory presented in~\cite{fajp}. We consider differential systems with
rapidly varying coefficients for which the origin is an equilibrium
point of neutral type. We state and prove a result concerning the
existence of the asymptotic phase for small solutions. We also state
and prove a result which is analogous to the Pliss reduction principle
for autonomous differential systems~\cite{plis}. We follow the
approach of Palmer~\cite{palm2} to these questions, though we will see
that some preparatory work is necessary in order to apply his
arguments.\par

It will be clear that the methods which we apply can be used to state
and prove results concerning the existence of the asymptotic phase and
the validity of the Pliss reduction principle for other classes of
non-autonomous dynamical systems. In this regard, we can refer to the
contributions of Janglajew~\cite{jang}, P{\"o}tzsche~\cite{poet}, and
Reinfelds-Janglajew~\cite{reja}.\par

This paper is organized as follows. In Section~\ref{sec_attractors},
we elaborate and compare the notions of attractor mentioned above. In
Section~\ref{sec_asymp}, we consider the questions of the existence of
the asymptotic phase and the validity of the Pliss reduction principle
for differential systems with rapidly varying coefficients.

\section{Attractors}\label{sec_attractors}

As stated in the Introduction, we will work in the \emph{skew-product}
framework for non-autonomous differential systems. We now describe
that framework.\par

Let $\fP$ be a compact metric space, and let $\{\tau_t:t\in\R\}$ be a
flow on $\fP$. That is, for each $p\in\fP$, the map
$\tau_t:\fP\rightarrow\fP$ is a homeomorphism and, moreover, the
following conditions are satisfied:
\begin{enumerate}
\item[(i)] $\tau_0(p)=p$ for all $p\in\fP$;
\item[(ii)] $\tau_t\circ\tau_s=\tau_{t+s}$ for all $t,s\in\R$;
\item[(iii)] the map $\tau:\fP\times\R\rightarrow\fP$,
  $(p,t)\mapsto\tau_t(p)$, is continuous.
\end{enumerate}

We introduce some terminology related to the flow concept. If
$p\in\fP$, then the \emph{orbit} through $p$ is
$\{\tau_t(p):t\in\R\}$. The \emph{positive} (resp. \emph{negative})
\emph{semiorbit} through $p$ is $\{\tau_t(p):t\geq 0\,\mbox{(resp. }t\leq
0\mbox{)}\}$. A subset $A\subset\fP$ is said to be \emph{invariant} if
$\tau_t(A)\subset A$ for all $t\in\R$; it is said to be
\emph{positively} (resp. \emph{negatively}) \emph{invariant} if
$\tau_t(A)\subset A$ for all $t\geq 0$ (resp. $t\leq 0$). If
$p\in\fP$, the \emph{omega-limit set}
$\omega(p)=\{p_1\in\fP:\mbox{there exists a sequence
}t_n\to\infty\mbox{ such that }\tau_{t_n}(p)\to p_1\}$.\par
Next, let $d\geq 1$ be an integer. We define the concept of
\emph{skew-product local flow} on $\fP\times\R^d$. By this, we mean a
pair $(Y,\tilde\tau)$ which satisfies the conditions given below. We
write $Z=\fP\times\R^d$, and $z=(p,x)$ for a generic point of $Z$.
\begin{enumerate}
\item[(i)] $Y$ is an open subset of $Z\times\R$ which contains the set
  $\{(z,0):z\in Z\}$, such that, if $(z,t)\in Y$, then $(z,s)\in Y$
  for all $0\leq s\leq t$ if $t\geq 0$, and for all $t\leq s\leq 0$ if
  $t\leq 0$. The map $\tilde\tau:Y\to Z$ is continuous.
\item[(ii)] $\tilde\tau(z,0)=z$ for all $z\in Z$.
\item[(iii)] $\tilde\tau(\tilde\tau(z,s),t)=\tilde\tau(z,t+s)$
  whenever $\tilde\tau(z,s)$ and $\tilde\tau(\tilde\tau(z,s),t)$ are
  defined.
\item[(iv)] If $z=(p,x)\in Z$ and if $\tilde\tau(z,t)$ is defined,
  then $\tilde\tau(z,t)=(\tau(p),x_t)$ where $x_t\in\R^d$.
\item[(v)] The pair $(Y,\tilde\tau)$ is maximal with respect to the
  properties (i)-(iv).
\end{enumerate}

Condition (iv) is called the \emph{skew-product} property of
$\tilde\tau$. We will sometimes write
$\tilde\tau(z,t)=\tilde\tau_t(z)$ when $(z,t)\in Y$.\par

One can define a skew-product local flow beginning with an appropriate
family of non-autonomous differential systems. Let $|{\cdot}|$ denote
the Euclidean norm on $\R^d$. Let $f:\fP\times\R^d\to\R^d$ be a
continuous map such that, for each compact set $K\subset\R^d$, the
uniform Lipschitz constant
\[
\mbox{Lip}_K=\sup_{p\in\fP}\sup_{x_1\neq
  x_2}\left\{\frac{|f(p,x_1)-f(p,x_2)|}{|x_1-x_2|}\right\}
\]
is finite. Consider the family of differential systems
\addtocounter{equation}{1}
\begin{equation}\tag*{(\arabic{section}.\arabic{equation})$_p$}\label{eq_hull}
  x'=f(\tau_t(p),x)
\end{equation}
where $p$ ranges over $\fP$. For each $p\in\fP$ and each $x_0\in\R$,
let $x(t)$ be the solution of~\ref{eq_hull} such that $x(0)=x_0$. Set
$z_0=(p,x_0)\in Z$, and put $\tilde\tau(z_0,t)=(\tau_t(p),x(t))$ where
$x({\cdot})$ is the maximal solution of~\ref{eq_hull} with initial
value $x_0$. Let $Y=\{(z,t):\tilde\tau\mbox{ is well-defined}\}\subset
Z\times\R$. Then $(Y,\tilde\tau)$ is a skew-product local flow on
$\fP\times\R^d$.\par

It is well-known that, if $\tilde f:\R\times\R^d\to\R^d$ is uniformly
continuous on $\R\times K$ for each compact set $K\subset\R^d$, and if
$\tilde f$ is uniformly Lipschitz continuous in $x\in K$ for each such
set $\R\times K$, then $\tilde f$ gives rise to a compact metric space
$\fP$, a flow $(\fP,\{\tau_t\})$, and a function
$f:\fP\times\R^d\to\R^d$ satisfying the above conditions. See
e.g.~\cite{sell}.\par

Let $\fP$ be a compact metric space, and let $(\fP,\{\tau_t\})$ be a
flow. Let $\tilde\tau$ be a skew-product local flow on
$\fP\times\R^d$, and let $\pi:\fP\times\R^d\to\fP$ be the projection.

\begin{definition}
  A compact set $A\subset Z=\fP\times\R^d$ is said to be a
  \emph{Lyapunov attractor} if the following conditions are satisfied.
  \begin{enumerate}
  \item[(i)] If $z=(p,x)\in A$, then $\tilde\tau(z,t)$ is defined for
    all $t\in\R$.
  \item[(ii)] There is an open set $W\subset Z$ containing $A$ such
    that, if $z\in W$,then $\tilde\tau_t(z)$ is defined for all $t\geq
    0$, and the omega-limit set $\omega(z)\subset A$.
  \item[(iii)] If $V$ is an open set in $Z$ containing $A$, then there
    is an open set $V_1\subset V$, which contains $A$, such that, if
    $z\in V_1$, then $\tilde\tau_t(z)$ is defined and lies in $V$ for
    all $t\geq 0$.
  \item[(iv)] The projection $\pi(A)=\fP$.
  \end{enumerate}
\end{definition}

The last condition is imposed for reasons of convenience. In practice,
it will not usually entail any loss of generality, since one has the
option of replacing $\fP$ by $\pi(A)$.\par

Our first goal is to show that a compact subset
$A\subset\fP\times\R^d$ is a Lyapunov attractor if and only if it is a
\emph{pullback attractor} (\cite{carc}, \cite{chks}, \cite{klro},
\cite{sche}, \cite{schm}). We give a definition of the latter concept
which seems appropriate in the context of skew-product local flows. We
make use of the concept of Hausdorff distance~\cite{klro}. Let
$(Z,d_Z)$ be a compact metric space, and let $K_1,K_2$ be two nonempty
compact subsets of $Z$. Set $H_*(K_1,K_2)=\sup_{z_1\in
  K_1}\inf\{d_Z(z_1,z_2):z_2\in K_2\}$, then set
\[
H(K_1,K_2)=\max\{H_*(K_1,K_2),H_*(K_2,K_1)\}.
\]
The quantity $H(K_1,K_2)$ is the \emph{Hausdorff distance} between
$K_1$ and $K_2$. If $z\in Z$ and if $K\subset Z$ is compact, we abuse
notation slightly and write $H(z,K)$ in place of $H(\{z\},K)$.\par

Let us introduce the following notation: if $Z=\fP\times\R^d$ and
$X\subset Z$, then $X_p=X(p)=(\{p\}\times\R^d)\cap X$. Further, let
$d_{\fP}$ be a metric on $\fP$ compatible with its topology, and
define the metric $d$ on $Z=\fP\times\R^d$ by
$d(z_1,z_2)=d_{\fP}(p_1,p_2)+|x_1-x_2|$ whenever $z_1=(p_1,x_1)$ and
$z_2=(p_2,x_2)$.

\begin{definition}\label{def_pull}
  Let $A\subset Z=\fP\times\R^d$ be a compact invariant set such that
  $\pi(A)=\fP$. Say that $A$ is a \emph{pullback attractor} if there
  is an open subset $W\subset Z$, which contains $A$, with the
  following property: let $D\subset W$ be a compact set such that
  $D_p\neq\emptyset$ for each $p\in\fP$; then
  \[
  \limsup_{t\to\infty}\sup_{p\in\fP}H(\tilde\tau_t(D_{\tau_{-t}(p)}),A)=0.
  \]
\end{definition}

Our definition is related to that given in \cite{sche},
\cite{schm}. However, our ``universe'' of sets $\{D\}$ is chosen so as
to take account of the skew-product framework in which we work.\par

The following result is quite natural; however, so far as we know it
has not appeared in the literature.

\begin{proposition}\label{lya_iff_pull}
  Let $A\subset\fP\times\R^d$ be a compact invariant set. Then $A$ is
  a Lyapunov attractor if and only if it is a pullback attractor.
\end{proposition}
\begin{proof}
  Suppose first that $A$ is a Lyapunov attractor. Let $W\subset
  Z=\fP\times\R^d$ be an open neighborhood of $A$ such that
  $\omega(z)\subset A$ for all $z=(p,x)\in W$. We claim that the
  pullback condition is satisfied by the family $\{D\}$ of compact
  subsets $D$ of $W$ such that $D_p\neq\emptyset$, $p\in\fP$.\par
  To see this, let $D\subset W$ be a compact subset such that
  $D_p\neq\emptyset$, $p\in\fP$. Suppose for contradiction that there
  exist a number $\varepsilon>0$ and sequences
  $t_n\to\infty,\,p_n\in\fP$ such that
  \[
  H(\tilde\tau_{t_n}(D(\tau_{-t_n}(p_n))),A)\geq\varepsilon,\,n\in\N.
  \]
  This means that, for each $n\geq 1$, there exists $x_n\in
  D(\tau_{-t_n}(p_n))$ such that, if
  $z_n=\tilde\tau_{t_n}(\tau_{-t_n}(p_n),x_n)$, then
  \begin{equation}\tag*{(*)}\label{pos_haus}
    H(z_n,A)\geq\varepsilon.
  \end{equation}
  On the other hand, let $\delta>0$ be a number such that, if $z\in Z$
  and $H(z,A)\leq\delta$, then $H(\tilde\tau_t(z),A)\leq\varepsilon/2$
  for all $t\geq 0$. Such a number exists because $A$ is a Lyapunov
  attractor. If $z_*\in D$, then $\omega(z_*)\subset A$. Hence there
  exists a time $t_*=t_*(z_*)>0$ such that, if $t\geq t_*$, then
  $H(\tilde\tau_t(z_*),A)<\delta$. There is a neighborhood
  $N_*=N_*(z_*)$ of $z_*$ in $W$ such that, if $z_1\in N_*$, then
  $H(\tilde\tau_t(z_1),A)<\delta$, and by compactness of $D$ we can
  find a fixed time $t_f$ such that, if $z\in D$, then for some time
  $t_z\in(0,t_f]$ there holds $H(\tilde\tau_{t_z}(z),A)<\delta$. But
  this implies that, if $t\geq t_f$, then
  $H(\tilde\tau_t(z),A)\leq\varepsilon/2$. This is inconsistent with
  the condition~\ref{pos_haus}. So indeed $A$ is a pullback
  attractor.\par

  Now let us suppose that the pullback condition is satisfied by
  $A$. Let $W\subset\fP\times\R^d$ be a neighborhood of $A$ for which
  the pullback condition holds. We claim that, if $z=(p,x)\in W$, then
  $\omega(z)\subset A$. For suppose not. Then there exist
  $\varepsilon>0$ and a sequence $t_n\to\infty$ such that
  \begin{equation}\tag*{(**)}\label{haus_tau_pos}
    H(\tilde\tau_{t_n}(z),A)\geq\varepsilon.
  \end{equation}
  Write $z_n=(p_n,x_n)=\tilde\tau_{t_n}(z)$, so that
  $\tilde\tau_{-t_n}(z_n)=z$. Then $x\in D_p$ for some compact set
  $D\subset W$ such that $D_{p_1}\neq\emptyset$ for all $p_1\in\fP$
  (one can take, for example, $D=A\cup\{z\}$). But
  then~\ref{haus_tau_pos} contradicts the pullback condition, so
  indeed $\omega(z)\subset A$.\par
  We now claim that the Lyapunov stability condition holds. For
  suppose not; then there is a number $\varepsilon>0$ with the
  following property: there are sequences
  $\{\delta_n\}\subset(0,\infty)$, $\{z_n\}\subset Z$, and
  $\{t_n\}\subset(0,\infty)$ such that $\delta_n\to 0$,
  $H(z_n,A)\leq\delta_n$, and
  \begin{equation}\tag*{(***)}\label{haus_zn_pos}
    H(\tilde\tau_{t_n}(z_n),A)\geq\varepsilon.
  \end{equation}
  There is no loss of generality in assuming that
  $\tilde\tau_{t_n}(z_n)\in W$.\par
  Now, however, $z_n=\tilde\tau_{-t_n}(\tilde\tau_{t_n}(z_n))$. Fix a
  compact neighborhood $D$ of $A$ in $W$. Then $z_n$ lies in $D$ for
  all large $n$. But then~\ref{haus_zn_pos} violates the pullback
  condition. This shows that $A$ is indeed a Lyapunov attractor, and
  completes the proof of Proposition~\ref{lya_iff_pull}.
\end{proof}

Next we consider the following situation. Let $V\subset\R^d$ be an
open set. Let $\mathcal C_V=\{c:\fP\to V:c\mbox{ is bounded and
  continuous}\}$ with the usual metric
\[
\rho(c_1,c_2)=\sup_{p\in\fP}\{|c_1(p)-c_2(p)|\}.
\]
Then $(\mathcal C_V,\rho)$ is a complete metric space. Suppose that,
for some $t_0>0$, there holds $\tilde\tau_{t_0}(\fP\times
V)\subset\fP\times V$. In this case, one can define a map
$\xi:\mathcal C_V\to\mathcal C_V$ as follows:
\begin{equation}\label{map_xi}
  (p,\xi(c)(p))=\tilde\tau_{t_0}(\tau_{-t_0}(p),c(\tau_{-t_0}(p))).
\end{equation}
We will suppose that the map $\xi$ is a \emph{contraction} with
respect to $\rho$. That is, there exists a constant $\alpha<1$ with
the property that
\[
\rho(\xi(c_1),\xi(c_2))\leq\alpha\rho(c_1,c_2)
\]
for all $c_1,c_2\in\mathcal C_V$.\par

Under these conditions, $\xi$ admits a unique fixed point
$a\in\mathcal C_V$: $\xi(a)=a$. Let
\[
A=\{(p,a(p)):p\in\fP\}.
\]
Then $A\subset\fP\times V\subset Z$ is a compact set, and it is
natural to conjecture that $A$ is a Lyapunov attractor for
$\tilde\tau$. We will verify that this is indeed the case.

\begin{proposition}\label{fix_then_lya}
  Let $V\subset\R^d$ be an open set, and suppose that $t_0>0$ is a
  number such that $\tilde\tau_{t_0}(\fP\times V)\subset\fP\times
  V$. Suppose that the map $\xi:\mathcal C_V\to\mathcal C_V$ defined
  by~{\upshape(\ref{map_xi})} is a contraction in $\mathcal C_V$ with
  contraction constant $\alpha<1$. Then the fixed point $a\in\mathcal
  C_V$ of $\xi$ has the property that
  \[
  A=\{(p,a(p)):p\in\fP\}
  \]
  is a Lyapunov attractor for $\tilde\tau$.
\end{proposition}
\begin{proof}
  We first show that $A$ is $\tilde\tau$-invariant. This is not
  immediately clear because $\tilde\tau$ is only assumed to be a local
  flow. Let $U\subset\R^d$ be an open set such that, if
  $z=(p,x)\in\fP\times V$, then $\tilde\tau_t(z)\in\fP\times U$ for
  all $t\geq 0$. Define $\mathcal C_U$ in the same way in which
  $\mathcal C_V$ was defined; then $\mathcal C_V\subset\mathcal
  C_U$. Fix a pair of integers $n\geq 1$, $0\leq r\leq n$, and set
  $u=rt_0/n$. Define a map $\xi_1:\mathcal C_V\to\mathcal C_U$ via the
  relation
  \[
  (p,\xi_1(c)(p))=\tilde\tau_{u}(\tau_{-u}(p),c(\tau_{-u}(p))),\,c\in\mathcal
  C_V.
  \]
  Then the iterate $\xi_1^{(k)}$ is defined for each $k\geq 1$. Note
  that $\rho(\xi_1(a),a)=\rho(\xi_1(\xi_1^{(n)}(a)),\xi_1^{(n)}(a))$
  $\leq\alpha^r\rho(\xi_1(a),a)$, so $a$ is a fixed point of
  $\xi_1$. This means that $\tilde\tau_u$ leaves $A$ invariant for all
  $u$ as above. By continuity of $\tilde\tau$, one has that
  $\tilde\tau_t(A)\subset A$ for all $t\in [0,t_0]$, and it follows
  that $A$ is positively $\tilde\tau$-invariant.\par
  We can now globally invert the local flow $\tilde\tau$ on $A$ by
  setting
  $\tilde{\tilde\tau}_{-t}(p,a(p))=(\tau_{-t}(p),a(\tau_{-t}(p)))$,
  $t\geq 0,p\in\fP$. We omit the proof that $\tilde{\tilde\tau}$
  coincides with $\tilde\tau$ on $A$.\par
  Let us now show that $A$ is a pullback attractor. In fact we will
  show that $W=\fP\times V$ satisfies the condition of
  Definition~\ref{def_pull}. For this, let $D\subset W$ be a compact
  set such that $D_p\neq\emptyset$ for each $p\in\fP$. We must show
  that
  \[
  \limsup_{t\to\infty}\sup_{p\in\fP}H(\tilde\tau_t(D(\tau_{-t}(p))),A)=0.
  \]
  Suppose for contradiction that there exists $\varepsilon>0$ with the
  following property: there are sequences $t_n\to\infty$, $p_n\in\fP$,
  and $x_n\in D(\tau_{-t_n}(p_n))$ such that
  \begin{equation}\label{haus_pn_pos}
    H(\tilde\tau_{t_n}(\tau_{-t_n}(p_n),x_n),A)\geq\varepsilon.
  \end{equation}
  For each $n\in\N$, there is a continuous map $c_n:\fP\to V$ such
  that $c_n(\tau_{-t_n}(p_n))=x_n$. For each $n$, let $k_n$ be the
  largest integer such that $k_nt_0\leq t_n$. We have
  \[
  \rho(\xi^{(k_n)}(c_n),a)=\rho(\xi^{(k_n)}(c_n),\xi^{(k_n)}(a))\leq\alpha^{k_n}\rho(c_n,a),
  \]
  and the right-hand side tends to zero as $n\to\infty$. Taking
  account of the definition of $\xi$ and of the continuity of
  $\tilde\tau$, we see that~(\ref{haus_pn_pos}) is violated for
  sufficiently large $n$. This completes the proof of
  Proposition~\ref{fix_then_lya}.
\end{proof}

\section{Asymptotic phase and reduction principle}\label{sec_asymp}

In this section, we state and prove the principle of asymptotic phase
and the Pliss reduction principle in a form appropriate for ordinary
differential systems with rapidly varying coefficients. We first
discuss the integral manifold theory for such a system as it is
presented in~\cite{fajp}.\par

Let $I\subset\R$ be a compact interval containing $\varepsilon=0$ in
its interior. Let $f:\R\times\R^d\times I\to\R^d$ be a function with
regularity and recurrence properties which will be specified
presently. Consider the non-autonomous differential system
\begin{equation}\label{rapid_eq}
  x'=|\varepsilon|f(t,x,\varepsilon),\;t\in\R,\,x\in\R^d,\,\varepsilon\in I
\end{equation}
where $\varepsilon$ is small. We write $|\varepsilon|$ instead of
$\varepsilon$ before $f$ because one might want, for instance, to
study the change of stability of an equilibrium point
of~(\ref{rapid_eq}) when $\varepsilon$ passes through zero.\par

The integral manifold theory for~(\ref{rapid_eq}) can be formulated in
an elegant way when $f$ depends on $t$ in a ``uniquely ergodic''
manner. We pause to explain what this means.\par
Let $\fP$ be a compact metric space and let $\{\tau_t\}$ be a flow on
$\fP$. We review some basic notions of ergodic theory. Let $\mu$ be a
regular Borel probability measure on $\fP$. Then $\mu$ is said to be
\emph{invariant} if, for each Borel set $B\subset\fP$ and each
$t\in\R$, there holds $\mu(\tau_t(B))=\mu(B)$. An invariant measure
$\mu$ is said to be \emph{ergodic} if it satisfies the following
indecomposibility condition: for each Borel set $B\subset\fP$, the
property
\[
\mu(\tau_t(B)\triangle B)=0\mbox{ for all }t\in\R\mbox{
  (}\triangle=\mbox{symmetric difference)}
\]
implies that either $\mu(B)=0$ or $\mu(B)=1$. It is known~\cite{nest}
that there exists at least one regular Borel measure $\mu$ on $\fP$
which is ergodic (with respect to $\{\tau_t\}$).\par

Let us now impose the following hypotheses, which will be valid in all
that follows. First, we assume that the flow $(\fP,\{\tau_t\})$ admits
a \emph{unique} invariant measure $\mu$, which is then necessarily
ergodic. Second, we assume that there exist: (i) a continuous map
$\tilde f:\fP\times\R^d\times I\to\R^d$ and (ii) a point $\tilde
p\in\fP$ such that
\[
f(t,x,\varepsilon)=\tilde f(\tau_t(\tilde p),x,\varepsilon)
\]
for all $t\in\R,\,x\in\R^d$, and $\varepsilon\in I$. This means in
effect that $f({\cdot},x,\varepsilon)$ has recurrence properties which
are reflected in the structure of the flow $(\fP,\{\tau_t\})$. Third,
we assume that $\tilde f(p,0,\varepsilon)=0$ for all
$p\in\fP,\,\varepsilon\in I$.\par

There are well-known conditions on $f$ which ensure the existence of
objects $\fP,\{\tau_t\},\mu$, and $\tilde f$ which satisfy the above
conditions. For example, if $f$ is almost periodic in $t$, uniformly
on each set of the form $K\times I$ where $K\subset\R^d$ is compact,
then the above conditions are fulfilled. See~\cite{fajp}
(also~\cite{sell} and many others references) for a discussion of this
matter. We will say no more about it; rather, we let
$\fP,\{\tau_t\},\mu$, and $\tilde f$ be as above, then write $f$
instead of $\tilde f$, and consider the family of differential systems
\addtocounter{equation}{1}
\begin{equation}\tag*{(\arabic{section}.\arabic{equation})$_p$}\label{rapid_hull}
  x'=|\varepsilon|f(\tau_t(p),x,\varepsilon)
\end{equation}
which includes the $\varepsilon$-dependent
equation~(\ref{rapid_eq}). We generally will not indicate explicitly
the dependence of the family~\ref{rapid_hull} on $\varepsilon\in
I$.\par

We proceed to outline the integral manifold theory for a ``uniquely
ergodic family''~\ref{rapid_hull} as it is worked out
in~\cite{fajp}. We will assume throughout that, for some $r\geq 1$,
the function $x\mapsto f(p,x,\varepsilon):$ $\R^d\to\R^d$ is of class
$C^r$ for each $p\in\fP$ and $\varepsilon\in I$. We further require
that the Fr{\'e}chet derivatives $\partial_x^kf:\fP\times\R^d\times
I\to L^k(\R^d,\R^d)$ are continuous, $0\leq k\leq r$. Here
$L^k(\R^d,\R^d)$ is the usual space of $\R^d$-valued, $k$-linear maps
defined on $\R^d\times\cdots\times\R^d$ ($k$ times).\par

We introduce the average
\[
\bar f(x,\varepsilon)=\int_\fP f(p,x,\varepsilon)d\mu(p)=
\lim_{|t|\to\infty}\frac{1}{t}\int_0^tf(\tau_s(p),x,\varepsilon)ds.
\]
Because of the unique ergodicity of $(\fP,\{\tau_t\})$, the limit on
the right is uniform on $\fP\times K\times I$ for each compact subset
$K\subset\R^d$. Note that $\bar f(0,\varepsilon)=0$ for all
$\varepsilon\in I$. Let us write
\begin{equation}\label{l_plus_n}
  \bar f(x,\varepsilon)=\bar l_\varepsilon x+\bar n_\varepsilon(x)
\end{equation}
where $\bar n_\varepsilon(x)=o(|x|)$ as $x\to 0$, uniformly in
$\varepsilon\in I$. Set $\varepsilon=0$; we impose a condition on the
matrix $\bar l_0$.
\begin{hypothesis}
  The set of eigenvalues of $\bar l_0$ is the union of two disjoint
  nonempty subsets, namely
  \[
  \begin{split}
    \Sigma^{(0)}&=\{\lambda\in\C:\lambda\mbox{ is an eigenvalue of
    }\bar
    l_0\mbox{ and Re}(\lambda)=0\},\\
    \Sigma^{(-)}&=\{\lambda\in\C:\lambda\mbox{ is an eigenvalue of
    }\bar l_0\mbox{ and Re}(\lambda)<0\}.
  \end{split}
  \]
\end{hypothesis}
Let $\beta>0$ be a number such that Re$(\lambda)<-\beta$ for all
$\lambda\in\Sigma^{(-)}$. Let $L^{(-)}\subset\R^d$ be the intersection
with $\R^d$ of the sum of the generalized eigenspaces of $\bar l_0$
which correspond to eigenvalues in $\Sigma^{(-)}$. Define $L^{(0)}$ in
the analogous way. Let $Q_0:\R^d\to\R^d$ be the projection with image
$L^{(-)}$ and kernel $L^{(0)}$.\par
Next write
\[
f(p,x,\varepsilon)=l_\varepsilon(p)x+n_\varepsilon(p,x)
\]
for $p\in\fP,x\in\R^d,\varepsilon\in I$. Consider the family of linear
systems \addtocounter{equation}{1}
\begin{equation}\tag*{(\arabic{section}.\arabic{equation})$_p$}\label{lin_hull}
  x'=|\varepsilon|l_\varepsilon(\tau_t(p))x.
\end{equation}
If $\varepsilon\neq 0$, the change of variables $s=|\varepsilon|t$
transforms~\ref{lin_hull} into \addtocounter{equation}{1}
\begin{equation}\tag*{(\arabic{section}.\arabic{equation})$_p$}\label{lin_hull2}
  \frac{dx}{ds}=l_\varepsilon(\tau_{s/|\varepsilon|}(p))x.
\end{equation}
Moreover,~\ref{rapid_hull} is transformed into
\addtocounter{equation}{1}
\begin{equation}\tag*{(\arabic{section}.\arabic{equation})$_p$}\label{rapid_hull2}
  \frac{dx}{ds}=l_\varepsilon(\tau_{s/|\varepsilon|}(p))x+n_\varepsilon(\tau_{s/|\varepsilon|}(p),x).
\end{equation}

We recall a basic definition~\cite{copp}, \cite{sase}.
\begin{definition}\label{def_ed}
  The family of equations~\ref{lin_hull2} is said to have an
  \emph{exponential dichotomy} (\emph{ED}) over $\fP$ if the following
  conditions hold. Let $\Phi_p(s)$ be the fundamental matrix solution
  of~\ref{lin_hull2}; then there are constants $k>0,\sigma>0$ together
  with a continuous, projection-valued function $p\mapsto Q_p:$
  $\R^d\to\R^d$, $Q_p^2=Q_p$, such that the following estimates hold:
  \[
  \begin{split}
    |\Phi_p(u)Q_p\Phi_p(s)^{-1}|&\leq ke^{-\sigma(u-s)},\,u\geq s,\\
    |\Phi_p(u)(I-Q_p)\Phi_p(s)^{-1}|&\leq ke^{\sigma(u-s)},\,u\leq s.
  \end{split}
  \]
\end{definition}

We introduce the \emph{dynamical spectrum} of the
family~\ref{lin_hull2}. For this, let $\lambda\in\R$ and consider the
translated family \addtocounter{equation}{1}
\begin{equation}\tag*{(\arabic{section}.\arabic{equation})$_p$}\label{trans_hull}
  \frac{dx}{ds}=[-\lambda I+l_\varepsilon(\tau_{s/|\varepsilon|}(p))]x.
\end{equation}
Then the dynamical spectrum $\Sigma(\varepsilon)$ of the
family~\ref{lin_hull2} is
\[
\Sigma(\varepsilon)=\{\lambda\in\R:\mbox{the family~\ref{trans_hull}
  does \emph{not} admit ED over }\fP\}.
\]

We can use a basic perturbation theorem of Sacker and Sell~\cite{sase}
to prove the following result.
\begin{theorem}\label{th_spectrum}
  Let $\beta$ be as above, and let $\alpha>0$ be a real number. There
  exists $\varepsilon_1>0$ such that, if
  $0<|\varepsilon|\leq\varepsilon_1$, then
  $\Sigma(\varepsilon)=\Sigma^{(0)}(\varepsilon)\cup\Sigma^{(-)}(\varepsilon)$
  where
  $\Sigma^{(0)}(\varepsilon)\subset\{\lambda\in\R:|\lambda|<\alpha\}$
  and
  $\Sigma^{(-)}(\varepsilon)\subset\{\lambda\in\R:\lambda<-\beta\}$.
\end{theorem}

In the following developments, we will assume that $\alpha<\beta$ and
that $I\subset[-\varepsilon_1,\varepsilon_1]$. We then have that
$\Sigma^{(0)}(\varepsilon)\cap\Sigma^{(-)}(\varepsilon)=\emptyset$ for
all $\varepsilon\in I$.\par
Next, fix $\lambda_*\in(-\beta,-\alpha)$, so that the family of
translated equations
\[
\frac{dx}{ds}=[-\lambda_*I+l_\varepsilon(\tau_{s/|\varepsilon|}(p))]x
\]
admits an ED over $\fP$. This family is of course~\ref{trans_hull}
with $\lambda=\lambda_*$. Let us write $Q_{p,\varepsilon}$ for the
dichotomy projection of this family of equations; see
Definition~\ref{def_ed}. It turns out that, if $\varepsilon\in I$,
then $Q_{p,\varepsilon}$ does not depend on the choice of
$\lambda_*\in(-\beta,-\alpha)$~\cite{sase}. Furthermore, one has the
following important continuity result; see~\cite{copp}, \cite{sase}.
\begin{theorem}\label{q_cont}
  The mapping
  \[
  (p,\varepsilon)\mapsto\left\{
    \begin{array}{ll}
      Q_{p,\varepsilon}&p\in\fP,\,0\neq\varepsilon\in I,\\
      Q_0&p\in\fP,\,\varepsilon=0,
    \end{array}
  \right.
  \]
  is continuous on $\fP\times I$.
\end{theorem}
Let us write $Q_{p,0}=Q_0$ for all $p\in\fP$.\par

We can now describe the integral manifold theory for the nonlinear
family~\ref{rapid_hull2}. Let us write
$L^{(0)}(p,\varepsilon)=\mbox{Ker}(Q_{p,\varepsilon})\subset\R^d$ and
$L^{(-)}(p,\varepsilon)=\mbox{Im}(Q_{p,\varepsilon})\subset\R^d$. Let
$e=\dim L^{(0)}$ so that $d-e=\dim L^{(-)}$. It follows from
Theorem~\ref{q_cont} that $\dim L^{(0)}(p,\varepsilon)=e$ and $\dim
L^{(-)}(p,\varepsilon)=d-e$ for all $(p,\varepsilon)\in\fP\times
I$. If $\delta>0$ is a real number, let
$B_\delta=\{x\in\R^d:|x|\leq\delta\}$, then set
$L^{(0)}(\delta,p,\varepsilon)=L^{(0)}(p,\varepsilon)\cap B_\delta$,
$L^{(-)}(\delta,p,\varepsilon)=L^{(-)}(p,\varepsilon)\cap B_\delta$
for each $(p,\varepsilon)\in\fP\times I$. Let us note that, for each
$0\neq\varepsilon\in I$, the family of nonlinear
equations~\ref{rapid_hull2} induces a skew-product local flow on
$\fP\times\R^d$. Precisely, if $p\in\fP$ and $x_0\in\R^d$, set
$\tilde\tau^\varepsilon_2(p,x_0,s)=\varphi(s)$ where $\varphi(s)$ is the
maximal solution of~\ref{rapid_hull2} satisfying
$\varphi(0)=x_0$. Then $\tilde\tau^\varepsilon$ is a local
skew-product flow on $\fP\times\R^d$ which covers the ``fast flow''
$\tau^\varepsilon$ on $\fP$ defined by
$\tau^\varepsilon_s(p)=\tau_{s/|\varepsilon|}(p)$. Say that a subset
$Y\subset\fP\times\R^d$ is \emph{locally invariant} with respect to
$\tilde\tau^\varepsilon$ if, for each $y=(p,x)\in Y$, there exists
$s_0>0$ such that, if $|s|<s_0$, then $\tilde\tau^\varepsilon(y,s)\in
Y$.

\begin{proposition}[\cite{fajp}]\label{int_man}
  There exist positive numbers $\varepsilon_2\leq\varepsilon_1$ and
  $\delta\in\R$ together with a family of maps
  \[
  h_{p,\varepsilon}:L^{(0)}(\delta,p,\varepsilon)\longrightarrow
  L^{(-)}(\delta,p,\varepsilon),\,p\in\fP,\,0<|\varepsilon|\leq
  \varepsilon_2
  \]
  such that the following conditions hold.
  \begin{enumerate}
  \item[(a)] If $0\neq|\varepsilon|\leq\varepsilon_2$, then the set
    \[
    M_\varepsilon=\bigcup_{p\in\fP}\{(p,x):x=u+h_{p,\varepsilon}(u),\,u\in
    L^{(0)}(\delta,p,\varepsilon)\}
    \]
    is a locally invariant subset of $\fP\times\R^d$ with respect to
    the local flow $\tilde\tau^\varepsilon$. We call $M_\varepsilon$
    an \emph{integral manifold} of the
    family~{\upshape\ref{rapid_hull2}}.
  \item[(b)] Introduce the ``center bundle''
    $E_\delta^{(0)}=\{(p,u,\varepsilon):u\in
    L^{(0)}(\delta,p,\varepsilon),p\in\fP,0\neq|\varepsilon|\leq\varepsilon_2\}$
    and the ``stable bundle'' $E_\delta^{(-)}=\{(p,u,\varepsilon):u\in
    L^{(-)}(\delta,p,\varepsilon),p\in\fP,0\neq|\varepsilon|\leq\varepsilon_2\}$. Then
    the map $(p,u,\varepsilon)\mapsto (p,h_{p,\varepsilon}(u),\varepsilon):$
    $E_\delta^{(0)}\to E_\delta^{(-)}$ is continuous. Moreover, for
    each $p\in\fP$ and $0\neq|\varepsilon|\leq\varepsilon_2$, the map
    $u\mapsto h_{p,\varepsilon}(u):$
    $L^{(0)}(\delta,p,\varepsilon)\to\R^d$ is of class $C^r$. For each
    $k=0,1,\ldots,r$, the Fr{\'e}chet derivatives
    $\partial^k_uh_{p,\varepsilon}$,
    $p\in\fP,\,0\neq|\varepsilon|\leq\varepsilon_2$, are
    well-defined. The collection of maps $\{h_{p,\varepsilon}\}$ is
    $F^r$ in the sense of Foster~\cite{fost} for each $\varepsilon$
    with $0\neq|\varepsilon|\leq\varepsilon_2$. This means that, for
    each $k=1,\ldots,r$ one has the following statement. If
    $(p_n,u_n,\eps_n)\to(p,u,\eps)$ in $E^{(0)}_\delta$, and if
    $x_n^{(1)}\to x^{(1)},\ldots,x_n^{(k)}\to x^{(k)}$ are convergent
    sequences in $\R^d$, then
    $\partial^k_{u_n}h_{p_n,\eps_n}(x_n^{(1)},\ldots,x_n^{(k)})\to \partial^k_{u}h_{p,\eps}(x^{(1)},\ldots,x^{(k)})$.
  \item[(c)] For each $p\in\fP$ and each
    $0\neq|\varepsilon|\leq\varepsilon_2$, one has
    $h_{p,\varepsilon}(0)=\partial_xh_{p,\varepsilon}(0)=0$. Thus the
    manifold $M_\varepsilon\cap(\{p\}\times\R^d)$ is tangent to
    $L^{(0)}(p,\varepsilon)$ at the origin $x=0$.
  \end{enumerate}
\end{proposition}

One actually has that the family $\{h_{p,\varepsilon}\}$ extends in a
$C^r$-way to $\varepsilon=0$. It is fairly clear what this means; in
any case see~\cite{fajp}.\par

From now on, we assume that
$I\subset[-\varepsilon_2,\varepsilon_2]$. Our next goal is the
following. Fix $\varepsilon\in I$ with $\varepsilon\neq 0$. We want to
introduce (time-dependent) coordinates $u\in\R^e,\,v\in\R^{d-e}$ in
$\R^d$ such that, for each $p\in\fP$, the family of linear
systems~\ref{lin_hull2} assumes the block-diagonal form
$\frac{du}{ds}=au,\,\frac{dv}{ds}=bv$. This is because we wish to
apply Palmer's methods \cite{palm2} to the study of the asymptotic
phase and the Pliss reduction principle relative to the integral
manifold $M_\varepsilon$. It turns out, however, that in general we
cannot achieve a block-diagonalization \emph{over $\fP$}. Instead we
must introduce an appropriate ``extension'' of the flow
$(\fP,\{\tau_s^\varepsilon\})$ to do so. We now discuss this question
in more detail.\par

It is convenient to assume that the original flow $(\fP,\{\tau_t\})$
is \emph{minimal}, i.e. that the orbit $\{\tau_t(p):t\in\R\}$ is dense
in $\fP$ for each $p\in\fP$. (This implies that, for each
$\varepsilon\neq 0$ in $I$, the fast flow
$(\fP,\{\tau_s^\varepsilon\})$ is minimal.) See \cite{elli} for a
study of minimal flows. The condition of minimality entails little
loss of generality for the following reason. The flow
$(\fP,\{\tau_t\})$ admits by assumption just one invariant measure
$\mu$, which is therefore ergodic. Let $\fP_\mu$ be the topological
support of $\mu$; that is, $\fP_\mu$ is the complement in $\fP$ of the
union over all open sets $O\subset\fP$ satisfying $\mu(O)=0$. Then
$\fP_\mu$ is invariant and $(\fP_\mu,\{\tau_t\})$ is minimal.\par

We assume form now on that $(\fP,\{\tau_t\})$ is minimal and uniquely
ergodic (one says that the flow is \emph{strictly ergodic}). This
implies that $(\fP,\{\tau_s^\varepsilon\})$ is strictly ergodic for
each $0\neq\varepsilon\in I$.\par

An \emph{extension} $(\fQ,\{T_t\})$ of $(\fP,\{\tau_t\})$ consists of
a compact metric space $\fQ$, a flow $\{T_t\}$ on $\fQ$, and a
continuous surjective map $\pi:\fQ\to\fP$ such that
\[
\tau_t(\pi(q))=\pi(T_t(q)),\,t\in\R,\,q\in\fQ.
\]
One says that $\pi$ is a \emph{homomorphism} of the flows
$(\fQ,\{T_t\})$ and $(\fP,\{\tau_t\})$. We will show that, if
$0\neq\eps$ is sufficiently small, then there is an extension of
$(\fP,\{\tau_s^\eps\})$ with respect to which one can
block-diagonalize the family~\ref{lin_hull2} after it has been
``lifted'' to $\fQ$ in an appropriate way. We discuss the appropriate
concept of lifting. To lighten the notation, we will write
$\{\tau_s\}$ instead of $\{\tau_s^\eps\}$ when $\eps\in I$ has been
fixed.\par

Suppose that $0\neq\eps\in I$, and that $(\fQ,\{T_s\})$ is an
extension of $(\fP,\{\tau_t\})$ with flow homomorphism
$\pi:\fQ\to\fP$. Set $l(q)=l_\eps(\pi(q))$ and
$n(q,x)=n_\eps(\pi(q),x)$. Consider the family of linear equations
\addtocounter{equation}{1}
\begin{equation}\tag*{(\arabic{section}.\arabic{equation})$_q$}\label{linear}
  x'=l(T_s(q))x
\end{equation}
and the family of nonlinear equations \addtocounter{equation}{1}
\begin{equation}\tag*{(\arabic{section}.\arabic{equation})$_q$}\label{nonlinear}
  x'=l(T_s(q))x+n(T_s(q),x).
\end{equation}
It is natural to view these families as lifts of the
families~\ref{lin_hull2} and~\ref{rapid_hull2}. Clearly the statements
of Theorems~\ref{th_spectrum}, \ref{q_cont} and \ref{int_man} can be
modified as to be valid for the lifted families~\ref{linear}
and~\ref{nonlinear}. We will not write out these modified versions of
Theorems~\ref{th_spectrum}-\ref{int_man}.\par

For each integer $n\geq 1$, let $\M_n$ be the set of $n\times n$ real
matrices with the usual operator norm. If $\varphi:\fQ\to\M_n$ is a
map, let $|\varphi|_0=\sup\{|\varphi(q)|:q\in\fQ\}$.

\begin{theorem}\label{th_change_block}
  There is a positive number $\eps_3\leq\eps_2$ such that, if
  $0\neq|\eps|\leq\eps_3$, then there exist a minimal extension
  $(\fQ,\{T_s\})$ of $(\fP,\{\tau_s\})$ together with a continuous
  function $\sigma:\fQ\to\M_d$ with the following properties. First,
%
%
%
%
%
%
%
%
%
%
%
%
%
%
%
  the map $\sigma^{-1}:\fQ\to\M_d$ is well-defined and
  continuous. Second, the map $\sigma':\fQ\to\M_d:$
  $\sigma'(q)=\left.\frac{d}{ds}\sigma(T_s(q))\right|_{s=0}$ is
  well-defined and continuous. Third, for each $q\in\fQ$ the change of
  variables
  \[
  x=\sigma(T_s(q))y
  \]
  transforms equation~\ref{linear} into the block form
  \addtocounter{equation}{1}
  \begin{equation}\tag*{(\arabic{section}.\arabic{equation})$_q$}\label{block_sys}
    y'=\left(
      \begin{array}{cc}
        a(T_s(q))&0\\
        0&b(T_s(q))
      \end{array}
    \right)y,
  \end{equation}
  where $a:\fQ\to\M_e$ and $b:\fQ\to\M_{d-e}$ are continuous. Fourth,
  $|a|_0\leq|l|_0$ and $|b|_0\leq|l|_0$.
\end{theorem}
\begin{proof}
  The first step involves an application of some results presented in
  Coppel's lecture notes (\cite[pp.  37-41]{copp}). See also
  Daletskii-Krein~\cite{dakr}.\par

  Let us note that, if $\lambda\in(-\beta,-\alpha)$ is fixed, and if
  $0\neq\eps\in I$, then the family~\ref{trans_hull} admits an ED over
  $\fP$ with a dichotomy constant $k$ which does not depend on
  $\eps$. Moreover, we can use Theorem~\ref{q_cont} to determine a
  positive number $\eps_3\leq\eps_2$ such that, if
  $0<|\eps|\leq\eps_3$, then for every $p\in\fP$, the
  angle~\cite{dakr} between $\mbox{Im}(Q_{p,\eps})$ and
  $\mbox{Im}(Q_0)$ is less than $\pi/6$, and the angle between
  $\mbox{Ker}(Q_{p,\eps})$ and $\mbox{Ker}(Q_0)$ is less than
  $\pi/6$.\par

  Taking account of these facts, we can apply Lemma 3 on p. 41
  of~\cite{copp} to determine a constant $\theta\geq 0$, which is
  independent of $p\in\fP$ and $0\neq|\eps|\leq\eps_3$, for which the
  following statements can be verified.\par

  Let $\lambda\in(-\beta,-\alpha)$, $\bar p\in\fP$, and
  $0\neq\bar\eps\in I$ be fixed. Let us write $l(s)=-\lambda
  I+l_{\bar\eps}(\tau_s^{\bar\eps}(\bar p))$. Consider the linear
  system
  \begin{equation}\label{lin_no_hull}
    \frac{dx}{ds}=l(s)x.
  \end{equation}
  There is a continuously differentiable function $\sigma:\R\to\M_d$,
  with continuously differentiable inverse $\sigma^{-1}$, with the
  following properties.
  \begin{enumerate}
  \item[(i)] The change of variables
    \[
    x=\sigma(s)y
    \]
    transforms equation~(\ref{lin_no_hull}) into the block-form
    \[
    y'=m(s)y,\,m(s)=\left(
      \begin{array}{cc}
        a(s)&0\\
        0&b(s)
      \end{array}
    \right),\,s\in\R
    \]
    where $a({\cdot})\in\M_e$ and $b({\cdot})\in\M_{d-e}$.
  \item[(ii)] $|m(s)|\leq|l(s)|$ for all $s\in\R$.
  \item[(iii)] One has $\sup_{s\in\R}|\sigma(s)|\leq\theta,\,
    \sup_{s\in\R}|\sigma(s)^{-1}|\leq\theta$.
  \item[(iv)]
    $\frac{d\sigma}{ds}=l(s)\sigma(s)-\sigma(s)l(s),\,s\in\R.$
  \end{enumerate}

  It follows from (iii) and (iv) that $\frac{d\sigma}{ds}$ is
  uniformly bounded, so $\sigma$ and $\sigma^{-1}$ are uniformly
  continuous functions. It then follows from (iv) that
  $\frac{d\sigma}{ds}$ is (bounded and) uniformly continuous.\par

  The second step involves a Bebutov-type flow and basic methods of
  topological dynamics. We view $\sigma$ as an element of the space
  $\cC=\{c:\R\to\M_d:c\mbox{ is bounded and continuous}\}$. We equip
  $\cC$ with the topology of uniform convergence on compact
  sets. There is a flow $\{T_s\}$ on $\cC$ induced by the
  translations: thus $T_s(c)(u)=c(s+u)$ for $c\in\cC$ and
  $s,\,u\in\R$. This is a flow of Bebutov-type~\cite{bebu}. Since
  $\sigma$ is uniformly continuous, the closure
  $\cC_\sigma=\mbox{cls}\{T_t(\sigma):t\in\R\}$ is compact; clearly
  $\cC_\sigma$ is $\{T_t\}$-invariant.\par

  If $c\in\cC_\sigma$, set $C(c)=c(0)$. Then $C:\cC_\sigma\to\M_d$ is
  continuous. One can verify that $C^{-1}:\cC_\sigma\to\M_d:$
  $C^{-1}(c)=c(0)^{-1}$ is well-defined and continuous. Moreover, a
  basic result of Analysis can be used to show that the map $C'$
  defined by $C'(c)=\frac{dc}{ds}(0)$ is well-defined and
  continuous.\par

  Let $p\in\fP$ be a general point, and let $\{t_n\}\subset\R$ be a
  sequence such that $\tau_{t_n}(\bar p)\to p$. There is a subsequence
  $\{t_k\}\subset\{t_n\}$ such that $T_{t_k}(\sigma)$ converges in
  $\cC$, say to $\sigma_p$. Let $l_p(s)=-\lambda
  I+l_{\bar\eps}(\tau_s^{\bar\eps}(p))$. One can check that, if $l_p$
  is substituted for $l$ and $\sigma_p$ is substituted for $\sigma$,
  then all statements (i)-(iv) above are valid when $m(s)$ is
  substituted by
  \[
  m_p(s)=\sigma_p(s)^{-1}l_p(s)\sigma_p(s)-\sigma_p(s)^{-1}\frac{d\sigma_p}{ds}(s).
  \]
  Caution: the function $m_p$ is not \emph{uniquely} determined by
  $p$, since different subsequences $\{t_k\}$ may give rise to
  different limit functions $\sigma_p$.\par

  Now let $\fQ\subset\cC_\sigma$ be an invariant set such that
  $(\fQ,\{T_s\})$ is minimal (see e.g.~\cite{elli}). There is a
  natural projection $\pi:\fQ\to\M_d:$ $\sigma(q)=C(q)$. By stepping
  through the definition, one can now verify all the assertions of
  Theorem~\ref{th_change_block}. This completes the proof.
\end{proof}

From now on, we suppose that $I\subset[-\eps_3,\eps_3]$. We note an
important corollary of the proof of Theorem~\ref{th_change_block}.
\begin{corollary}\label{block_bound}
  Let $0\neq\eps\in I$, and let $(\fQ,\{T_s\})$ and $\sigma$ be as in
  the statement of Theorem~\ref{th_change_block}.
  \begin{enumerate}
  \item[(i)] There is a constant $\theta\geq 0$, which does not depend
    on $q\in\fQ$ and $0\neq\eps\in I$, such that
    $|\sigma(q)|\leq\theta$ and $|\sigma^{-1}(q)|\leq\theta$.
  \item[(ii)] Let $\Phi_q^{(a)}(s)$ be the fundamental matrix solution
    of the system
    \[
    \frac{du}{ds}=a(T_s(q))u,
    \]
    and let $\Phi_q^{(b)}(s)$ be the fundamental matrix solution of
    the system
    \[
    \frac{dv}{ds}=b(T_s(q))v.
    \]
    There is a constant $k_1>0$, which does not depend on $q\in\fQ$
    and $0\neq\eps\in I$, such that
    \[
    \begin{split}
      |\Phi_q^{(a)}(s)|&\leq k_1e^{\alpha|s|},\,s\in\R,\\
      |\Phi_q^{(b)}(s)|&\leq k_1e^{-\beta s},\,s\geq 0.
    \end{split}
    \]
  \end{enumerate}
\end{corollary}
The number $k_1$ is only distantly related to the dichotomy constant
$k$ of the family~\ref{trans_hull} introduced in the proof of
Theorem~\ref{th_change_block}, and may be much larger. Nevertheless,
we will indicate $k_1$ by $k$ in the succeeding developments.\par

Next let $0\neq\eps\in I$. Introduce a minimal flow $(\fQ,\{T_s\})$
which satisfies the conditions of Theorem~\ref{th_change_block}. Let
$\sigma,a$, and $b$ be the functions of that theorem. Consider the
nonlinear family~\ref{nonlinear} obtained by lifting the
family~\ref{rapid_hull2} to $\fQ$. For each $q\in\fQ$, introduce the
change of variables $x=\sigma(T_s(q))y$ in
equation~\ref{nonlinear}. Set
$y=\left(\!\!\begin{array}{c}u\\v\end{array}\!\!\right)$ where
$u\in\R^e,\,v\in\R^{d-e}$, and $\R^d=\R^e\oplus\R^{d-e}$. Further set
\[
g(q,u,v,\eps)=\sigma^{-1}(q)n(q,\sigma(q)y),
\]
then put $g=\left(\!\!\begin{array}{c}g_1\\g_2\end{array}\!\!\right)$
where $g_1\in\R^e$, $g_2\in\R^{d-e}$. Then equations~\ref{nonlinear}
take the form \addtocounter{equation}{1}
\begin{equation}\tag*{(\arabic{section}.\arabic{equation})$_q$}\label{hull_block}
  \begin{split}
    \frac{du}{ds}&=a(T_s(q))u+g_1(T_s(q),u,v,\eps),\\
    \frac{dv}{ds}&=b(T_s(q))v+g_2(T_s(q),u,v,\eps).
  \end{split}
\end{equation}
Fix $\eps\in I$. Return to the spaces $L^{(0)}(p,\eps)$,
$L^{(-)}(p,\eps)$, $L^{(0)}(\delta,p,\eps)$, and
$L^{(-)}(\delta,p,\eps)$ of Proposition~\ref{int_man}. These spaces
can be lifted to $\fQ$ by setting
$L^{(0)}(q,\eps)=L^{(0)}(\pi(q),\eps)$, etc.; we have made an obvious
abuse of notation. The functions $h_{p,\eps}$ lift naturally to
functions $h_{q,\eps}:L^{(0)}(\delta,q,\eps)\to
L^{(-)}(\delta,q,\eps)$ where again we have abused notation. For a
given number $\delta>0$, let $\hat\delta=\delta\theta^{-1}$ where
$\theta$ is the constant of Corollary~\ref{block_bound}. Let us write
\[
\R^e_{\hat\delta}=\R^e\cap\{y\in\R^d:|y|\leq\hat\delta\}.
\]
We see that $\sigma(q)\R^e\subset L^{(0)}(q,\eps)$,
$\sigma(q)\R^{d-e}\subset L^{(-)}(q,\eps)$, and
$\sigma(q)\R^e_{\hat\delta}\subset L^{(0)}(\delta,q,\eps)$ for all
$q\in\fQ$.\par

Note that the solutions of the family~\ref{hull_block} of
non-autonomous differential systems determine a local skew-product
flow $\tilde T^\eps$ on $\fQ\times\R^d$.\par
Consider the functions
\[
\hat
h_{q,\eps}=\sigma(q)^{-1}h_{q,\eps}\sigma(q):\R^e_{\hat\delta}\longrightarrow\R^{d-e},
\]
and set
\[
\hat M_\eps=\bigcup_{q\in\fQ}\{(q,y):y=u+\hat
h_{q,\eps}(u),\,u\in\R^e_{\hat\delta}\}.
\]
Then $\hat M_\eps$ is a locally invariant subset of $\fQ\times\R^d$
with respect to the local flow $\tilde T^\eps$. Moreover, the
functions $\hat h_{q,\eps}$ satisfy conditions analogous to those
stated in parts (b) and (c) of Proposition~\ref{int_man}. In particular,
the collection of maps $\{\hat h_{q,\eps}:q\in\fQ\}$ is $F^r$ in the
sense of Foster, and $\hat h_{q,\eps}(0)=\partial_u\hat
h_{q,\eps}(0)=0$. We abuse notation still again, and write $\delta$
for $\hat\delta$, $h_{q,\eps}$ for $\hat h_{q,\eps}$, and $M_\eps$ for
$\hat M_\eps$.\par

We will now discuss results concerning the asymptotic phase and the
Pliss reduction principle which make reference to the set $M_\eps$
defined by the family~\ref{hull_block}. It will be clear that these
results imply corresponding statements which make reference to the set
$M_\eps$ as originally defined for the family~\ref{lin_hull2}.\par

The preceding constructions allow us to follow Palmer's arguments
\cite{palm2} in a fairly straightforward way. He works in the context
of autonomous differential systems, however. For the reader's
convenience, we sketch how his statements and proofs can be modified
so as to apply to our family of non-autonomous
equations~\ref{hull_block}.\par

Choose $\delta>0$ so that $|\partial_uh_{q,\eps}|\leq 1$ for all
$q\in\fQ$, $|u|\leq\delta$, $0\neq\eps\in I$. Then
$|h_{q,\eps}(u)|\leq|u|$ whenever $|u|\leq\delta$. Fix $0\neq\eps\in
I$. For each positive number $\Delta\leq\delta$, let $\omega(\Delta)$
be the maximum, as $q$ ranges over $\fQ$, of the upper bounds of the
norms of the derivatives $\partial_ug_1$, $\partial_vg_1$,
$\partial_ug_2$, $\partial_vg_2$, $\partial_uh_{q,\eps}$ in the set
$|u|\leq\Delta$, $|v|\leq\Delta$. Then $\omega(\Delta)$ decreases to
zero as $\Delta\to 0^+$.\par

If $q\in\fQ$, consider the linearization of~\ref{hull_block}:
\addtocounter{equation}{1}
\begin{equation}\tag*{(\arabic{section}.\arabic{equation})$_q$}\label{lin_block}
  \left(\!\!\begin{array}{c}\frac{du}{ds}\\\frac{dv}{ds}\end{array}\!\!\right)= 
  \left(\begin{array}{cc}a(T_s(q))&0\\0&b(T_s(q))\end{array}\right)
  \left(\!\!\begin{array}{c}u\\v\end{array}\!\!\right).
\end{equation}
It admits the fundamental matrix solution
\[
\left(\begin{array}{cc}\Phi_q^{(a)}(s)&0\\0&\Phi_q^{(b)}(s)\end{array}\right).
\]
By Corollary~\ref{block_bound}, there is a constant $k>0$, which does
not depend on $q\in\fQ$ and $0\neq\eps\in I$, such that
\[
\begin{split}
  |\Phi_q^{(a)}(s)|&\leq ke^{\alpha|s|},\,s\in\R,\\
  |\Phi_q^{(b)}(s)|&\leq ke^{-\beta s},\,s\geq 0.
\end{split}
\]

The conditions above are analogous to those stated on \cite[p. 274]{palm2}. Caution: 
our numbers $\alpha$ resp. $\beta$ play the roles of his $\beta$ resp. $\alpha$.

\begin{lemma}\label{lemma_approach}
  Let $0\neq\eps\in I,\,q\in\fQ$. Write $h(s,u)=h_{T_s(q),\eps}(u)$,
  $a(s)=a(T_s(q))$, $b(s)=b(T_s(q))$, $g_i(s,u,v)=g_i(T_s(q),u,v)$,
  $1\leq i\leq 2$, $s\in\R$, $u\in\R^e$, $v\in\R^{d-e}$. Let $\gamma$
  be a real number such that $0<\gamma<\beta-\alpha$.\par
  There is a positive number $\Delta\leq\delta/2$, which does not
  depend on the choice of $\eps$ and $q$, such that the following
  statements are valid. Let $S>0$, and let
  $\left(\!\!\begin{array}{c}u(s)\\v(s)\end{array}\!\!\right)$ be a
  solution of~{\rm\ref{hull_block}} which is defined on the interval
  $[0,S]$. Then the solution $\tilde u(s)$ of the system
  \begin{equation}\label{first_block}
    \frac{du}{ds}=a(s)u+g_1(s,u,h(s,u))
  \end{equation}
  with $\tilde u(s)=u(s)$ is defined on $[0,S]$. Moreover, $|\tilde
  u(s)|\leq 2\Delta$ and
  \[
  |u(s)-\tilde u(s)|+|v(s)-h(s,\tilde u(s))|\leq
  2k|v(0)-h(0,u(0))|e^{-(\beta-\gamma)s}
  \]
  for all $0\leq s\leq S$.
\end{lemma}
\begin{proof}
  Let us first show that $h(s,u)$ is a $C^1$ function of $(s,u)$. This
  is not immediately obvious, because Proposition~\ref{int_man} states
  only that $h$ is continuous in $s$. However, we will see that the
  local invariance of $M_\eps$ actually implies that $h$ is $C^1$ in
  $(s,u)$. Fix $\Delta<\delta$, and let
  $\R^e_\Delta=\{u\in\R^e:|u|\leq\Delta\}$.\par

  We will show that the partial derivatives $\partial_sh$,
  $\partial_uh$ exist at each point $(s_0,u_0)\in\R\times\R^e_\Delta$,
  and are jointly continuous on $\R\times\R^e_\Delta$. First note
  that, by Proposition~\ref{int_man}, the partial Fr{\'e}chet
  derivative $\partial_uh(s_0,u_0)$ is defined and continuous as
  $(s_0,u_0)$ varies over $\R\times\R^e_\Delta$. Moreover one has
  \[
  h(s_0,u_0+u)-h(s_0,u_0)-\partial_uh(s_0,u_0)u=o(u)
  \]
  where $o(u)/|u|\to 0$ as $u\to 0$, uniformly in $s_0\in\R$.\par

  We show that the partial derivative $\partial_sh$ exists and is
  continuous on $\R\times\R^e_\Delta$. For this, fix $s_0\in\R$ and
  $u_0\in\R^e_\Delta$. Let $v_0=h(s_0,u_0)$, and let
  $\left(\!\!\begin{array}{c}\bar u({\cdot})\\ \bar
      v({\cdot})\end{array}\!\!\right)$ be the solution
  of~\ref{hull_block} satisfying $\bar u(s_0)=u_0$, $\bar
  v(s_0)=v_0$. We have
  \[
  h(s_0+s,u_0)-h(s_0,u_0)=h(s_0+s,u_0)-h(s_0+s,\bar u(s_0+s))
  +h(s_0+s,\bar u(s_0+s))-h(s_0,u_0).
  \]
%
%
%
%
%
%
%
%
%
%
%
%
  However $h(s_0+s,u_0)-h(s_0+s,\bar
  u(s_0+s))=-\partial_uh(s_0+s,u_0)(\bar u(s_0+s)-u_0)+o(\bar
  u(s_0+s)-u_0)$ and $h(s_0+s,\bar u(s_0+s))-h(s_0-u_0)=\frac{d\bar
    v}{ds}(s_0)+o(s)$ where we use the invariance of $M_\eps$ to
  obtain the second relation. Letting $s\to 0$ we obtain
  \begin{equation}\tag*{(*)}\label{partial_h}
    \begin{split}
      \partial_sh(s_0,u_0)=&-\partial_uh(s_0,u_0)[a(s_0)u_0+g_1(s_0,u_0,h(s_0,u_0))]\\
      &+[b(s_0)h(s_0,u_0)+g_2(s_0,u_0,h(s_0,u_0))].
    \end{split}
  \end{equation}
  The explicit expression for $\partial_sh$ in~\ref{partial_h} shows
  that it is continuous in its arguments.\par

  Now we follow the arguments of \cite{palm2}. Choose $\Delta>0$ such
  that $\Delta\leq\delta/2$ and
  $4k^2\omega(2\Delta)\leq\min\{2\gamma,\beta-\alpha-\gamma,4k^2\}$. Let
  $\left(\!\!\begin{array}{c}u(s)\\v(s)\end{array}\!\!\right)$ be the
  solution of~\ref{hull_block} referred to in the statement of the
  present lemma. Write $z(s)=v(s)-h(s,u(s))$. Then
  \[
  \begin{split}
    \frac{dz}{ds}=&b(s)z(s)+g_2(s,u(s),v(s))-g_2(s,u(s),h(s,u(s)))\\
    &-\partial_uh(s,u(s))[g_1(s,u(s),v(s))-g_1(s,u(s),h(s,u(s)))]
  \end{split}
  \]
  where we used~\ref{partial_h}. One completes the proof if the lemma
  by mimicking the estimates in \cite{palm2}; we omit the details.
\end{proof}

Fix numbers $\gamma\in(0,\beta-\alpha)$ and $\Delta>0$ which satisfy
the conditions of Lemma~\ref{lemma_approach}. We retain the notation
introduced in the statement of Lemma~\ref{lemma_approach}.

\begin{corollary}\label{sol_in_man}
  Fix $0\neq\eps\in I$ and $q\in\fQ$. Let
  $\left(\!\!\begin{array}{c}u(s)\\v(s)\end{array}\!\!\right)$ be a
  solution of~{\rm\ref{hull_block}} such that $|u(s)|\leq\Delta$,
  $|v(s)|\leq\Delta$ for all $s\in\R$. Then $v(s)=h(s,u(s))$ for all
  $s\in\R$, and hence
  $\left(\!\!\begin{array}{c}u(s)\\v(s)\end{array}\!\!\right)\in
  M_\eps$ for all $s\in\R$.
\end{corollary}

This corollary may be proved by adapting the reasoning of
\cite[Proposition 1]{palm2}. The next statement is proved by
appropriate modification of the arguments of \cite[Proposition
2]{palm2}.
\begin{theorem}[Asymptotic phase]
  Let $\left(\!\!\begin{array}{c}u(s)\\v(s)\end{array}\!\!\right)$ be
  a solution of~{\rm\ref{hull_block}} which is defined and satisfies
  $|u(s)|\leq\Delta$, $|v(s)|\leq\Delta$ for all $s\geq 0$. Then there
  exists a solution $u_\infty(s)$ of~{\rm(\ref{first_block})} such
  that
  \[
  |u(s)-u_\infty(s)|+|v(s)-h(s,u_\infty(s))|\leq
  2k|v(0)-h(0,u(0))|e^{-(\beta-\gamma)s}
  \]
  for $s\geq 0$.
\end{theorem}

Observe that
$\left(\!\!\begin{array}{c}u_\infty(s)\\h(s,u_\infty(s))\end{array}\!\!\right)\in
M_\eps$ for all $s\geq 0$. Thus the solution
$\left(\!\!\begin{array}{c}u(s)\\v(s)\end{array}\!\!\right)$
``tracks'' a positive semiorbit in $M_\eps$ as $s\to\infty$.\par

Finally we have
\begin{theorem}[Pliss Reduction Principle]\label{pliss}
  Let $A\subset\fQ\times\R^d$ be a compact, $\tilde T^\eps$-invariant
  set such that, if $(q,y)\in A$ and
  $y=\left(\!\!\begin{array}{c}u\\v\end{array}\!\!\right)$, then
  $|u|\leq\Delta/2$ and $|v|\leq\Delta/2$. Then $A\subset
  M_\eps$. Moreover, if $A$ is a Lyapunov attractor relative to
  $M_\eps$, then it is a Lyapunov attractor relative to
  $\fQ\times\R^d$.
\end{theorem}
\begin{proof}
  It follows immediately from Corollary~\ref{sol_in_man} that
  $A\subset M_\eps$.\par

  We indicate how the second statement can be proved. Introduce the
  family of equations \addtocounter{equation}{1}
  \begin{equation}\tag*{(\arabic{section}.\arabic{equation})$_q$}\label{local_eq}
%
%
%
%
%
%
%
%
%
%
%
    \frac{du}{ds}=a(T_s(q))u+g_1(T_s(q),u,h_{T_s(q),\eps}(u)).
  \end{equation}
  Let $U=\{u\in\R^e:|u|<\delta\}$. The family~\ref{local_eq} gives
  rise to a local flow on $\fQ\times U$. The mapping $i:\fQ\times U\to
  M_\eps:$ $(q,u)\mapsto(q,h_{q,\eps}(u))$ is a diffeomorphism onto
  its image, which maps trajectories of the local flow in $\fQ\times
  U$ onto trajectories of the local flow in $M_\eps$. Let
  $A_0\subset\fQ\times U$ be the compact invariant subset which is the
  preimage of $A$ with respect to this diffeomorphism.\par

  One now argues as in \cite[Proposition 3]{palm2} to show that, if
  $A_0$ is a Lyapunov attractor relative to (the local flow on)
  $\fQ\times U$, then it is a Lyapunov attractor relative to
  $\fQ\times\R^d$. This clearly implies that the second statement of
  Theorem~\ref{pliss} is valid.
\end{proof}

\end{document}